\documentclass[12pt,centertags,oneside]{amsart}

\usepackage{amsmath,amstext,amsthm,amscd,typearea,hyperref,stmaryrd,mathabx}
\usepackage{amssymb}
\usepackage{a4wide}
\usepackage[mathscr]{eucal}
\usepackage{mathrsfs}
\usepackage{typearea}
\usepackage{charter}
\usepackage{pdfsync}
\usepackage[a4paper,width=16.2cm,top=3cm,bottom=3cm,lines=50]{geometry}

\numberwithin{equation}{section}

\newtheorem{thmspec}{\relax}

\newcommand{\cali}[1]{\mathscr{#1}}

\newcommand{\loc}{{loc}}
\newcommand{\ddc}{{dd^c}}

\newcommand{\dbar}{{\overline\partial}}
\newcommand{\ddbar}{{\partial\overline\partial}}

\newcommand{\dR}{{\rm dR}}

\renewcommand{\Re}{{\rm Re}}

\newcommand{\Cc}{\cali{C}}

\newcommand{\Sc}{\cali{S}}

\newcommand{\C}{\mathbb{C}}
\newcommand{\D}{\mathbb{D}}

\renewcommand{\P}{\mathbb{P}}

\title{Pseudoconvex domains with smooth boundary in projective spaces}

\author{Nessim Sibony}
\address{ Universit\'e
Paris-Saclay, CNRS,Laboratoire de Math\'ematiques d'Orsay, 91405 Orsay, France}

\address{and Korea  Institute for Advanced Study (KIAS), 85 Hoegiro, Dongdaemun-gu,
Seoul 02455, Republic of Korea}
\email{Nessim.Sibony@math.u-psud.fr}

\date{September 27, 2019}

\begin{document}

\begin{abstract}
Given  a pseudoconvex domain $U$  with $\Cc^1$-boundary  in $\P^n,$ $n\geq 3,$  we show that if $H^{2n-2}_{\dR}(U)\not=0,$ then
there is a strictly  psh function  in a neighborhood of $\partial U.$  We also    solve the $\dbar$-equation in $X=\P^n\setminus U,$  for data in $\Cc^\infty_{(0,1)}(X).$

We discuss Levi-flat domains in surfaces. If $Z$ is a real algebraic hypersurface in $\P^2,$  (resp a real-analytic  hypersurface with a
point of strict pseudoconvexity), then there is a strictly psh function in a neighborhood of $Z.$ 
\end{abstract}

\maketitle

\noindent
{\bf Classification AMS 2010:} Primary: 32Q28, 32U10;  32U40; 32W05;  Secondary 37F75


 \section{Introduction} \label{intro}


 In this paper  we discuss the pluri-potential theory on a smooth hypersurface  $Z$ in $\P^n.$ This includes the question of existence of positive closed (resp. $\ddc$-closed) currents supported on
  $Z$ and also the question of the  existence of a strictly plurisubharmonic  function (psh) in a  neighborhood of $Z.$  We will sometimes need a pseudo-convexity hypothesis on 
 a component of the complement. We also give some results in the case where $Z$ is a closed set satisfying some geometric assumptions. 
 
 Recall that  a complex  manifold of dimension $n$  is  strongly $q$-complete if it admits a  smooth exhaustion function $\rho$  whose Levi form  at each point has at least
 $(n-q+1)$ strictly positive eigenvalues. The main result in that theory is the following Theorem.

\renewcommand{\thethmspec}{Theorem}
\begin{thmspec}\label{T:AG}{\rm (Andreotti-Grauert \cite{AG})}.
Let $U$ be  a  strongly $q$-complete  manifold. Then for every  coherent analytic sheaf $\Sc$ over $U,$ $H^k(U,\Sc)=0,$ for $k\geq q.$

 In particular, if $k\geq q, H^{n,k}(U, \C)=0$ 
 \end{thmspec}
Indeed, for a holomorphic  bundle $E,$  $H^{n,k}(U,E)=H^k(U,\Lambda^{n,0} U\otimes E).$

Our main result will use the  above theorem.

\renewcommand{\thethmspec}{Theorem 1.1}
\begin{thmspec}
 Let $U$ be a domain in $\P^n,$ $n\geq 3,$  with $\Cc^1$ boundary.  Assume $U$ is strongly $(n-2)$-complete (i.e. it admits a smooth exhaustion function  whose  Levi form
 has at least $3$  strictly positive  eigenvalues  at each point).  Assume  also that  $H_\dR^{2n-2}(U,\C)\not=0.$
 Then there is a  strictly psh function  near the boundary of $U.$
\end{thmspec}
As  a  consequence of Theorem 1.1 and of Proposition 2.1 below, we  get the  following result.

\renewcommand{\thethmspec}{Corollary 1. 2}
\begin{thmspec}
 In $\P^n,$ $n\geq 3,$   there is  no $\Cc^1$  hypersurface $Z$ such that the two components $U^\pm$ of $\P^n\setminus Z$  are both
   strongly $(n-2)$-complete and one them, say $U^-$,  satisfies $H^{2}_\dR(U^-)=0$.
\end{thmspec}

Observe that in $\P^n$ the hypothesis $H^{2}_\dR(U^-)=0$ for a domain with  $\Cc^1$ boundary implies that $H_\dR^{2n-2}(U^+,\C)\not=0.$ See the proof of Proposition 2. 1 below.

Y. T. Siu has proved the following result, \cite{Siu2}.

\renewcommand{\thethmspec}{Theorem}
\begin{thmspec}{\rm  (Siu).}
 In $\P^n,$ $n\geq 3,$  there is no  Levi-flat  hypersurface $Z$ i.e both sides are $1$-complete (exhaustion with $n$ strictly positive eigenvalues).
\end{thmspec}
The meaning of  Levi flat is that  the Levi form  of a defining function $r$ for $Z,$ is  identically   zero  on the complex tangent  space.  Since the Levi problem
 has  a positive  solution  in $\P^n,$ this  implies that both   components $ U^\pm$ are Stein  and hence strongly $1$-complete. In particular   they are strongly $(n-2)$-complete for $n\geq 3.$ 
 
 We also  study  the  solvability of the $\dbar$-equation (resp.$\ddbar$-equation) on a  pseudo-concave set $X$ in the $\Cc^\infty$ category, i.e. we assume that $\P^n\setminus X$ is Stein. We use the H\"ormander duality method  (see  \cite{BS} for example).

The existence of $\ddc$-closed currents on $Z,$ gives an obstruction to the resolution of the $\dbar$- equation in the smooth category, see Theorem 6.1.

 In the last section  we discuss the  same problem in $\P^2.$  So far  it is  not known if there are smooth Levi-flat hypersurfaces in $\P^2.$
 
 \noindent {\bf Acknowledgements.} It is a pleasure thank Bo Berndtsson and Tien-Cuong Dinh for
their insightful comments and the referee for his questions.

\section{Real hypersurfaces  in $\P^n$ }


Let $Z$ be a  real hypersurface in $\P^n,$ $n\geq 2,$ of class $\Cc^1.$ Since $\P^n$ is  simply connected, $\P^n\setminus Z$ has two  components $U^\pm.$
 Denote by $\omega$ a K{\"a}hler form of mass one
on $\P^n.$

\renewcommand{\thethmspec}{Proposition 2.1}
\begin{thmspec}
If $n=2$,  the form $\omega$ is $d$-exact either on a neighborhood of $\overline{U^-}$ or on a neighborhood of $\overline{ U^+}.$

 If $n\geq 3,$ then: $H_\dR^{2n-2}(U^{+},\C)\not=0,$ iff $H_\dR^2\overline{ (U^-)}=0$.  In particular  the form $\omega$ is  $d$-exact, in a  neighborhood of $\overline{ U^-}$.
 
\end{thmspec}
\proof
Assume first $n=2$.
If $\omega$  is  not $d$-exact in a neighborhood of  $\overline{U^+},$ there is a 2-cycle $\sigma^+$ in  $\overline{U^+}$  which is non-trivial. Using that $Z$ is $\Cc^1,$
we can retract $\sigma^+$as  a cycle in  ${U^+}.$ Similarly,  we   would get a nontrivial  cycle $\sigma^-$ in  ${U^-}.$ But  by Poincar\'e  duality, $\sigma^+\sim  c_+\omega$ 
and  $\sigma^-\sim  c_-\omega,$  with $c_+ ,c_-$ non zero. Since $\sigma^+$ and $\sigma^-$ are  disjoint, we get $\sigma^+\smile \sigma^-=0,$ and hence  $c_+ c_-=0,$ a contradiction. 

 If $n\geq 3$ and  $H_\dR^{2n-2}(U^{+},\C)\not=0.$ Then there is a $(2n-2)$- cycle $\sigma^+$ in $U^+$, with $\sigma^+\sim c_+\omega ^{2n-2}$.
  If $H_\dR^2 (\overline {U^-}) \not=0,$ 
 we would construct as above a non-trivial $2$-cycle in $U^-$ and get a contradiction as above.
\endproof

Observe that  we cannot have  $\omega^{n-1},$  $d$-exact near  $\overline{U^+}$ and $\omega,$ $d$-exact near $\overline{U^-}.$ Otherwise, assume that  $\omega^{n-1}=d(\phi_+)$ near $ \overline{U^+}$ and  $\omega=d(\phi_-)$ near
$ \overline{U^-}.$  Then, $\omega^{n-1}-d(\chi_+\phi_+)$ and 
 $\omega-d(\chi_-\phi_-)$  would have  disjoint support  for appropriate  cut-off functions $\chi_\pm,$ contradicting  that the cup-product  should be non-zero.
 
  For  a compact $X$ in $\P^n,$  we will write $\mathcal H^2(X)=0,$ if there is an open  neighborhood
$V\supset X,$ such that the de Rham cohomology group  $H^2_\dR(V,\C)=0.$ 
 
\renewcommand{\thethmspec}{Proposition 2.2}
\begin{thmspec}
 Assume that $\omega$ is  $d$-exact  in a neighborhood of $\overline{U^-}.$ 
 \begin{enumerate}
   \item  There is  no closed current $A$ of order zero and  dimension  $2p$    supported on  $\overline{U^-},$  such that  $\{A\}\not=0.$ Here $ \{.\}$ denotes the de Rham cohomology class. In particular, there is no non-zero positive closed current supported on $\overline{U^-}.$ 
   \item  If $T$ is  a real  $\ddbar$-closed current (non-closed) of bidegree $(1,1)$  supported on $\overline{U^-}$  and $\{T\}\not=0.$ Then    there is a $\dbar$-closed $(2,0)$ holomorphic form
   non-identically  zero  in $U^+.$
  \end{enumerate}
 \end{thmspec}
\proof

 Since $\omega=d\varphi,$ near $\overline{U^-},$ we have
$$
\langle A,\omega^p\rangle=\langle A\wedge \omega^{p-1}, d\varphi  \rangle=0.
$$
Hence $\{A\}=0$.  If $A$ is positive closed non-zero, necessarily  $\{A\} \not=0$.

  It is also possible to define the cohomology class of a $\ddbar$-closed current in a compact K{\"a}hler manifold. It suffices to use the $\ddbar$- lemma, and Poincar\' e's duality. 
  
  When $T$ is a real $\ddbar$-closed current of bi-degree $(1,1),$ it follows from basic Hodge theory,  see \cite{FS} that when $\{T\}=\{\omega\},$ then
$$
T=\omega+\partial \sigma+\overline{\partial\sigma},
$$
where $\sigma$ is  a  $(0,1)$-current. Define
$$
T_c:= \omega+\partial\sigma +\overline{\partial\sigma}+\dbar \sigma+\partial\bar\sigma.
$$
It is  easy to  check that  $d T_c=0.$ If $\partial\bar \sigma=0$ on $\overline{U^+},$ then   $T_c$ is  closed  and supported on $\overline{U^-}.$
Since  $\{T_c\}=\{\omega\}$  we get a contradiction. It follows  that $\partial\bar \sigma$ is    not   identically  zero  in $\overline{U^+}.$
But  $\dbar \partial \bar\sigma=-\partial T$ is  supported  on $\overline{U^-}.$
Hence, $\dbar (\partial \bar\sigma)=0$ on ${U^+}.$ Therefore,  $\partial \bar\sigma$ is a holomorphic $(2,0)$-form.
\endproof

\renewcommand{\thethmspec}{Corollary 2.3}
\begin{thmspec}
Let $Z$ be a $\Cc^1$  hypersurface  in $\P^n.$ There is no positive closed current of bi- degree $(1,1)$ supported on $Z.$
\end{thmspec}
\proof
If $n=2$, this follows from Proposition 2.2, since we can consider that $Z$ bounds $\overline{U^-},$ and that $H_\dR^2(\overline{U^-})=0.$

Assume $n\geq3.$ Let $T$ a positive closed current of bidegree $(1,1)$ supported on $Z.$ Fix a point $p\notin Z$ and consider subspaces $L_p$ of  co-dimension$(n-2)$ through
 $p.$  For almost all  $L_p$, $L_p\cap Z$ is of class $\Cc^1,$ and the slice of $T$ at $L_p$ is a positive closed current.  Hence using the case $n=2,$ almost all slices vanish. This is true for all $p$ out of $Z.$ It follows then, from slicing theory, that $T=0.$
\endproof

\renewcommand{\thethmspec}{Remark 2.4}
\begin{thmspec} \rm

If $T$ is a closed and flat current, supported on $Z$ then  slicing theory is valid. We obtain that  $\{T\}=0.$ Indeed the class of a slice is the slice of the class.
\end{thmspec}

 \renewcommand{\thethmspec}{Corollary 2.5}
\begin{thmspec}
 Let $\overline{U^-},$ be a domain in $\P^2$ with $\Cc^1$ boundary. Assume, $H_\dR^2({U^-})=0.$ Then, there is a neighborhood of $\overline {U^-},$  which is Kobayashi hyperbolic.

\end{thmspec}
\proof
Otherwise, in an arbitrary neighborhood of $\overline{U^-},$ we will have a non-constant holomorphic image of 
$\C.$ This will permit to construct an Ahlfors current. In particular, we will have a positive closed current of mass one on $\overline{U^-}.$ Hence it's cohomology class is non-zero, contradicting Proposition 2.2.
\endproof

\section{Stricly psh functions near a compact $X$ and currents.}

 Let $(M,\omega)$ be  a  complex   Hermitian  manifold of dimension $n.$
 Let $X\Subset M$ be a compact set. We are interested in some general facts about the existence 
 of strictly psh functions near $X.$ Strictly psh functions are the starting point in order to use 
 H\"ormander's $L^2$ estimates, see for example J. J. Kohn \cite{Ko}. In particular, they permit to
 prove regularity at the boundary, for the $\dbar$-equation.
 
 \renewcommand{\thethmspec}{Proposition 3.1}
\begin{thmspec}
Let $X\Subset  M$ be  a compact set.  There is a positive  $\ddbar$-closed current $T$ of bi-dimension $(1,1)$ supported on $X$ iff  there is   no smooth  strictly psh  function   $u$ in a  neighborhood of  $X.$
Moreover, for any  $\Cc^2$ function $r,$
vanishing on $X$,
any such current $T,$ satisfies the following equations:

  \begin{equation}
    \label{e:(1)}
    \begin{split}
   T\wedge \partial r=0,\qquad T\wedge \ddbar r=0. 
    \end{split}
   \end{equation}

 \end{thmspec}
 \proof
 The proof is essentially the same as in \cite{S3} Proposition 2.1. There, $X=\partial U$ and $U$
 is smooth and pseudoconvex. Indeed, pseudoconvexity is not needed. The result is used  for arbitrary $X,$ in \cite{S3} Theorem 4.3. 
 
 If $u$  is  a  $\Cc^2$  psh  function in a  neighborhood  of $X$ and  $T$ is  a  positive  current 
supported on  $X,$ then

$$
\langle T,i\ddbar u  \rangle=\langle i\ddbar T,u   \rangle.
$$
So if $i\ddbar T=0,$ we get  that  $T=0,$ near every point where  $u$ is strictly psh. Hence if there is a strictly psh function near $X$, then $T=0$.
 
 We just show that any  positive $\ddbar$-closed current $T$ of mass one supported on $X$ satisfies the above relations. The proof of the other assertions is identical to the one  in \cite{S3}, mainly Hahn-Banach Theorem.
 
 Since $T$ is $\ddbar$-closed  then  
  $\langle T, i\ddbar r^2\rangle=0.$ Expanding and using that $T$ is positive and is supported on $\{r=0\},$we get: 
  $T\wedge i\partial r\wedge\overline{\partial} r=0.$
Therefore,  $T\wedge \partial r=0.$

Let $\chi$ be  a  smooth non-negative  function   with  compact  support. Using that
$T\wedge \partial r=0,$ we get that:
$$
0=\langle T,i\ddbar (\chi r)  \rangle=\langle T, \chi i\ddbar r   \rangle.
$$ 
Since $\chi$ is arbitrary, the measure,  $T\wedge i\ddbar r=0.$ 
\endproof

\renewcommand{\thethmspec}{ Remark 3.2}
\begin{thmspec} \rm

If $T$ is positive and $\ddbar$-closed on $M,$ then the calculus can be extended to continuous psh functions near $Support(T)$ \cite{DS3}. In fact if a function $u$ is continuous on  $Support(T)$  and is locally approximable, on  $Support(T),$  by continuous psh function, then
 $T\wedge i\ddbar u =0.$ 
 
  In particular, let $U$ be a pseudo-convex domain with boundary of class $\Cc^2$ in $\P^n.$ According to \cite {OS}, $U$ admits a bounded, strictly psh continuous exhaustion function $u.$
   Since $\partial U$ is of class $\Cc^2. $  For  $ p\in \partial U,$  $u$ is approximable by psh functions  in a fixed neighborhood of $p.$ Indeed, it suffices to push functions in the normal direction at $p.$ It follows that for $T$ positive of bi-dimension (1,1), $\ddbar$-closed  and supported on $\overline U,$ we have: $T\wedge i\ddbar u =0.$  Hence $T$ is supported on $\partial U.$
  \end{thmspec}

\renewcommand{\thethmspec}{Corollary 3.3}
\begin{thmspec}
Let $X\Subset  M$ be  a  compact subset. Assume there is no strictly psh function in a neighborhood of $X.$  Then there is  a compact $X_{\infty} \subset X$ with $X_{\infty} = \overline{\bigcup_{\alpha} X_{\alpha}},$ each $X_{\alpha},$ is compact connected and every continuous psh function,
in a neighborhood of $X_{\alpha}$ is constant on $X_{\alpha}.$ Moreover, there is a positive $\ddc$-closed current $T$ of bidimension $(1,1)$ such that $X_{\infty} =Support (T).$

 For any  compact $K \subset X$ with  $K \cap X_{\infty} = \emptyset,$ there are strictly psh functions near $K.$

\end{thmspec}
 \proof
  Assume there is no strictly psh function near $X.$ Then there is a positive  $\ddbar$-closed current $T$ of bi-dimension $(1,1)$ supported on $X,$  of mass $1.$ We can assume 
  $T$ is extremal and define $X'= Support (T).$ Then according to Proposition 4.2 in \cite {S3}, every continuous psh function near $X'$ is constant on $X'.$
  
  Let $\mathcal C_{1,1}$ denote the convex compact set of positive $\ddc$-closed currents supported on $X,$ of bi-dimension $(1,1)$ and of mass $1.$ Let $(T_{\alpha}),$ be the family of extremal elements in $\mathcal C_{1,1}.$ Let $X_{\alpha}: = Support (T_{\alpha}).$ Define 
$X_{\infty} = \overline{\bigcup_{\alpha} X_{\alpha}}.$  Since $(T_{\alpha})$ is extremal, it's support $X_{\alpha}$ is connected.

Let $(T_n), n \geq 1,$ be a dense sequence in $\mathcal C_{1,1}.$ Define : $T_{\infty}= \sum _{n}  2^{-n} T_n ,$ then $T_{\infty} \in \mathcal C_{1,1}$ and  $X_{\infty} =Support (T_{\infty}).$

As we have seen every $X_{\alpha}$ has the property that continuous psh function in a neighborhood of $X_{\alpha},$ is constant on  $X_{\alpha}.$
If $K \subset X$ and  $K  \cap X_{\infty} = \emptyset,$ then there are strictly psh functions near $K.$

 Indeed by Krein-Milman Theorem $\mathcal C_{1,1},$ is the closed convex hull of it's extremal elements. Hence there is no positive $\ddbar$-closed current $T$ of bi-dimension $(1,1)$ supported on $K.$  Proposition 3.1 implies that there is a strictly psh function near $K.$

\endproof
\renewcommand{\thethmspec}{Remarks 3.4}
\begin{thmspec} \rm
 \begin{enumerate}
 
\item Following \cite {Su}, one should call $X_{\infty} $ the  Poincar\'e set of $X$ and $ T_{\infty},$ a Poincar\'e current for $X.$

 There are many examples of the above decomposition in holomorphic  dynamics. It could happen that $X\setminus X_{\infty}, $ contains a biholomorphic image of $ \C^2,$  this is the case in the dynamics of H\'enon maps, if we take $X= \overline K^+,$ see \cite {DS4}.
 
  If $X$ is the the closure of the Torus in  Grauert's example for the Levi problem, as described in \cite {S3}, then there are uncountably many $X_{\alpha},$ each one being a real torus and also the closure of an image of $\C.$
  
In fact for a current $T_{\alpha}$ as above, a continuous $T_{\alpha}$ subharmonic functions (in the sense of \cite  {S1}) is necessarily constant. These are the function which are
 decreasing limits of $\Cc^2$ functions $u,$ satisfying $\ddc u \wedge T_{\alpha} \geq 0.$ This is a more intrinsic property, since it depends on the "complex directions" of  $T_{\alpha}$
 i.e. the infinitesimal complex structure of $X,$ independently of any smoothness assumption.
  \item 
A similar decomposition is given in \cite {S3} Theorem 4.3, for a domain $U$ admitting a continuous psh exhaustion function $\varphi.$ The obstruction to Steiness, is the existence
of positive $\ddbar$-closed currents supported on the level sets $\varphi = c.$ 

 The " Poincar\'e" decomposition of an arbitrary
domain $V \subset M,$ could be introduced following \cite {S3}, Definition 4.1. The Poincar\'e set  $V_{\infty}$ is the union of support of Liouville currents, i.e.  positive current of bi-dimension $(1,1),$ such that  $i\ddbar T=0$  and $T\wedge i\ddbar v=0,$
 for every  bounded continuous psh function in $V.$ 
 
 One can prove that $V_{\infty}$ is the support of a Liouville current $T_{\infty}$ and that $V_{\infty}$ is 1-pseudo-convex.

\item According to \cite {FS1} Corollary 2.6,  if $T$ a positive bi-dimension $(1,1)$, $\ddbar$-closed current, then $M \setminus Support(T)$ is $1$-pseudo-convex or with another terminology, $Support(T)$ is $1$-pseudo-concave. So $X_{\infty}$ is $1$-pseudo-concave. Hence the existence of a positive bi-dimension $(1,1),$ $\ddbar$-closed current, always implies the existence of a  $1$-pseudo-concave set. In fact there is a Poincar\'e decomposition $(X_{\alpha})$ of $Support(T)$ and each $X_{\alpha}$ is $1$-pseudo-concave.

In dimension $2$ and if $M$ is a surface where the Levi-problem has a positive solution, there is a smooth strictly psh  exhaustion function, on $M \setminus X_{\infty}.$ 
\item  
Consider a compact set  $X\Subset  M,$ and a closed pluripolar set $E \subset X.$ Assume $X\setminus E,$ satisfies the local maximum principle for continuous psh functions.
More precisely if $p \in X\setminus E,$ and $V_p$  is a neighborhood of $p,$ disjoint from $E.$ Then for any continuous psh function $u$ near $\overline V_p$ , we have
 $$ u(p) \leq \max_{z \in X\cap \partial V_p}  u(z).$$ 
 Then there is no strictly psh function on $X$ and hence there are extremal positive $\ddc$-closed currents supported on $X.$
 
 The proof is basically the same as in \cite {BS}.  Suppose there is a strictly psh function $u$ near $X.$ Assume it reaches it's maximum on $X$ at the point $p.$ 
 We can assume $p=0$ in a local chart, with local coordinates $z.$ Let $v$ be a psh  near $p,$ such that $v=-\infty$ on $E.$  For $\epsilon$
small enough, the function $u+ \epsilon v(z),$  will have a maximum at a  point $q \notin E$ near $p.$ So we can assume $0=p \notin E.$ Then for an appropriate cut-off function $\chi$, equal to $1$ near $0,$  the function  $u- \epsilon  \chi(z). \|z\|^2,$ will have a strict maximum at $p,$ contradicting the local maximum principle.

In particular if $X=E$ is pluripolar, either there is a strictly psh function near $E$ or it admits a Poincar\'e decomposition.
\item  Slodkowski  \cite {Sl} has shown, that a closed set $Z \subset M$ satisfies the local maximum principle iff it is $1$-pseudo-concave. 
\item
 It follows from Corollary 3.3 that smooth psh functions in a neighborhood of $X$ separate points in $X,$ iff there is a strictly psh function near $X.$ Indeed, if they separate points, ther is no $T_\alpha,$ hence there is a strictly psh function. For the converse, one can observe, that if $u$ is a strictly psh function near $X,$ then smooth function on a level set of $u$ can be extended to a strictly psh function.
\end{enumerate}

  \end{thmspec}

\renewcommand{\thethmspec}{Corollary 3.5}
\begin{thmspec}
Let $U\Subset  M$ be  a  domain  with $\Cc^2$ boundary.  Ler $r$ be a defining function for $\partial U$. Let $W\subset \partial U$ denote the set of points where the Levi-form is not positive or negative definite. Assume every component of $W$ is of 2-Hausdorff measure zero.  Then, there is a smooth  strictly psh  function  $u$ in a  neighborhood of  $\partial U.$

 \end{thmspec}
 \proof
 We can  assume that  
 $
 U:=\left\lbrace z\in U_1:\  r(z)<0  \right\rbrace
 $ where $U_1$ is a neighborhood of $\overline U.$ We have that $\partial r$ does not vanish on
 $\partial U.$
 
Assume there is no strictly psh function near $\partial U$. Let $T$ be an extremal positive $\ddbar$-closed current, of mass $1,$ supported on $\partial U$. 

 Recall that the Levi-form is defined on the complex tangent space of the boundary. If
$ \langle\partial r(z),t  \rangle=0,$ then the Levi-form at the point $z$ for the direction $t,$ is given by: $ \langle i\ddbar r(z),it\wedge\bar{t} \rangle.$ 

  At points of the boundary where $i\ddbar r>0$ or $i\ddbar r<0,$ on the complex tangent space, it follows  from the equation $T\wedge i\ddbar r=0,$ 
 that the current $T$ has no mass there, hence it is supported on $W.$ Since it is extremal it is supported on a component of $W.$  But as observed in \cite{BS}, positive $\ddbar$-closed currents give no mass to sets of 2-Hausdorff measure zero. So $T=0$ and the assertion follows.
\endproof

\renewcommand{\thethmspec}{Remark 3.6}
\begin{thmspec} \rm
Let  $X$ be a compact set in $M.$ Let $E$ be a closed subset of $X,$ of 2-Hausdorff measure zero. Assume that for every point $p\in X\setminus E$ there is a neighborhood $V_p$
of $p$ and a continuous psh function $u_p$ in $V_p$ peaking at $p$ on $X\cap V_p.$ Then there is a strictly psh function in a neighborhood of $X.$

Indeed, one can construct a continuous psh function, $v_p$ in a neighborhood of $X$, strictly psh at $p$ and peaking at $p$. Using Remark 3.2, one shows that a $\ddc$-closed current $T$ supported on $X$ has no mass near $p$. As above, it follows that $T=0.$
\end{thmspec}

A similar argument gives the following. Let $X$ be a real compact sub-manifold in $M$.  If the set $E$ of points in $X$ where there is a complex tangent is of 2-Hausdorff measure zero, then there is a smooth  strictly psh  function  $u$ in a  neighborhood of  $X.$ Indeed any positive $\ddbar$-closed  current of bi-dimension $(1,1),$ has to be supported on $E$.

\renewcommand{\thethmspec}{Corollary 3.7}
\begin{thmspec} Let $C$ be a compact connected real surface in $M.$ There is a strictly psh function near $C$ iff $C$ is not a complex curve.

\end{thmspec}
 \proof
 
 It is clear that if $C$ is a complex curve, there is no strictly psh function near $C$ (by maximum principle).
 
 Recall that the support of a positive 
 $\ddbar$-closed currents,  satisfies the local maximum principle for  local psh functions, \cite {S1} Theorem 3.2.

 Asssume $C$ is not a complex curve.  Let $T$ be a positive $\ddbar$-closed current of bi-dimension $(1,1)$ supported in $C.$  Since $C$ is a manifold, equations (3.1)  permit to consider $T$ as a current on $C.$ Let $E_c$ denote the set of points in $C$ where the tangent space is complex. Since $C$ is not a complex curve, then $E_c$ admits boundary points in $C.$  The current  $T$ is of bidimension $(1,1)$ and is supported on $E \subset E_c.$ Let $p$ be a boundary point of $E$ in $C.$ There are psh functions in a fixed neighborhood $V$ of $p$  with a unique peak point on $V$ near  
 $p.$  This contradicts the above local maximum principle. So there is no such $T.$
\endproof
\section{Constructing strictly psh functions}\label{S:Spsh}


We give a stronger version of Theorem 1.1. We do not assume, that $X$ is a domain with $\Cc^1$ boundary. When $X=\overline U{^-},$ is a domain with $\Cc^1$ boundary,
it is equivalent to assume that $ H_{dR}^2(U^-)=0$ or that $\mathcal H^2(X)=0,$ as explained in Proposition 2.1.
 
\renewcommand{\thethmspec}{Theorem 4.1}
\begin{thmspec}

Let $X$  be  a  compact set  in $\P^n,$ $n\geq 3,$ such that $\mathcal H^2(X)=0.$
Assume that the open set  $U^+:=\P^n\setminus X$ is   strongly $(n-2)$-complete. Then  there is  a  strictly psh   function  in a neighborhood of $X.$
 
\end{thmspec}
 
\proof
Let  $\Cc^\infty_{(0,1)}(X)$  denote  the space  of smooth forms of bidegree $(0,1)$ on $X.$
Here the smoothness is  in the Whitney sense  with  the  usual $\Cc^\infty$-Fr\'echet topology. Smooth functions on $X,$ in the Whitney sense, do extend as smooth 
functions in a neighborhood of $X.$ They admit also an intrinsic characterization using only the jet on $X$ i.e. the collection of derivatives.The jet extends and it is the jet of a smooth function. This permit to give the space a Fr\'echet topology, "uniform convergence on derivatives" \cite {MA}. 

The  dual space of  $\Cc^\infty_{(0,1)}(X)$ is the space of currents $R$ of bidegree $(n,n-1)$ on $\P^n,$ supported on $X.$

Let $\mathcal E$ denote the closure in  $\Cc^\infty_{(0,1)}(X)$ of $\{\dbar u\},$ $u\in\Cc^\infty(X).$  
We want to use the Hahn-Banach Theorem, to show that $\varphi^{0,1}$ is in $\mathcal E.$
Let $R$ be  a current, supported on $X,$ vanishing on the subspace $\mathcal E.$ We need to show that $\langle R,\varphi^{0,1}\rangle =0.$ 
Since the current $R$ is supported on $X,$ and vanishes on the subspace $\mathcal E,$ then $\dbar R=0$ on 
$\P^n.$
It follows from, $H^{n,n-1}(\P^n)=0,$ that there is  $S$ of bidegree $(n,n-2)$ such that  $R=\dbar S.$ Moreover, $S$ is smooth on $U^+.$ Indeed, $S$ is constructed using canonical solutions of the Hodge Laplacean, which satisfy the same regularity as the right hand side. Here  we are using the  local regularity in Hodge theory, \cite {DE}.

The  Andreotti-Grauert Theorem implies that on $U{^+},$ since $\dbar S=0,$ there is a form $B$ such that $S=\dbar B,$ on $U{^+}.$
Let $V\supset X$ be  an open  neighborhood of $X$   such that on $V,$ 
$$
\omega=d\varphi=\partial \varphi^{0,1}+\dbar \varphi^{1,0},\quad \dbar \varphi^{0,1}=0.
$$
Let $\chi$ be a cutoff function  with $\chi=1$ in a neighborhood of $U^+\setminus V$ and vanishing near $X.$ Then $R=\dbar S=\dbar[S-\dbar (\chi B)].$
Observe  that $S_1:= S-\dbar (\chi B)$ is  supported on $V,$   where $\dbar\varphi^{0,1}=0.$ Hence,
$$
\langle R,\varphi^{0,1}\rangle =\langle \dbar S_1, \varphi^{0,1}\rangle=-\langle S_1,\dbar \varphi^{0,1}\rangle=0.
$$
It  follows, by Hahn-Banach theorem, that $\varphi^{0,1}\in \mathcal E.$ Hence, there is a family $(u_\epsilon)$ of smooth functions  such that
$\dbar u_\epsilon\to \varphi^{0,1}$ in $\Cc^\infty(X).$ Then  $\omega=\lim_{\epsilon\to 0} i\ddbar \big({u_\epsilon-\bar u_\epsilon\over i}\big).$

As a consequence,  for $\epsilon>0$ small  enough and  $v_\epsilon:= {u_\epsilon-\bar u_\epsilon\over i},$   $i\ddbar v_\epsilon\geq  {1\over 2}\omega$ on $X.$
Hence, $v_\epsilon$ is  strictly psh  near $X.$
\endproof

To get Theorem 1.1, we should take $X=\overline U{^-}.$

  In order to prove Corollary 1.2, we will use the following version of the maximum principle, implicit in  \cite {S1}.
 \renewcommand{\thethmspec}{Lemma 4.2}
\begin{thmspec} 
Let $\rho$ be a function of class $\Cc^2$ in a neighborhood of a closed ball $ \overline B$ in $\C^k.$ Assume that for every point $ z\in B,$ there is a direction $t_z$ such that
$\langle i\ddbar \rho(z),it_z\wedge\bar{t_z} \rangle>0.$ Then there is no local maximum of $\rho$ in $B.$

\end{thmspec}

 \proof
 Assume by contradiction, that $\rho$ has a local maximum at a point $p\in B.$ Consider a complex disc $D_p$, at p in the direction $t_p.$ The restriction of $\rho$ to $D_p$ is strictly
 subharmonic on $D_p,$ near $p.$ It cannot have a local maximum at $p.$
 \endproof
 
 We now prove Corollary 1.2.
 
 \proof
 Assume to get a contradiction that the component $U^-,$ satisfies  $H^2(\overline{U^-})=0.$ According to Theorem 4.1, there is a  strictly psh  function  $v$ near $Z.$
 There is a point $z_0\in Z$ where $v(z_0)=\max_Z v.$ Using that $v$ is strictly psh,  we can assume the maximum at $z_0$ is  strict and $v(z_0)=0.$ Indeed, it suffices to add a negative small perturbation vanishing to second order at $z_0.$ So, $\left\lbrace v<0\right\rbrace,$ is a strictly  pseudoconvex  domain near $z_0.$ Hence there is a germ of complex hypersurface $W$ tangent to $U^-$ at $z_0$ and such that $W \setminus (z_0 )\subset U^+.$ Then $W \setminus( z_0) $
 is  strongly $(n-2)$-complete, with an exhaustion function $\rho,$  with $ 2$ strictly positive eigenvalues at each point, going to $+\infty $ at $z_0.$ Recall that the restriction of a 
 strongly $q$-complete function to a submanifold is still $q$-complete.
 
 Without loss of generality, we can assume that $W \setminus( z_0) $ is a pointed ball $B^{*}$ of dimension $(n-1).$  We can find a sequence of balls $(B_j)$ of dimension $(n-2)$ whose centers $(z_j)$ converge  to $z_0.$ On $B_j,$the function $\rho$ has a Levi-form, with one strictly positive eigenvalue. It satisfies the maximum principle, given in Lemma 4.2. Since $\rho,$ is uniformly
 bounded on $\bigcup (\partial B_j),$ it cannot converge to $+ \infty$ at $0.$  This finishes the proof . The last part shows that the pointed ball $B^{*}$ of dimension $(n-1),$ is not strongly $(n-2)$-complete.
 \endproof



\section{$\dbar$ equation on pseudo-concave sets}\label{S:dbar}


\renewcommand{\thethmspec}{Theorem 5.1}
\begin{thmspec}
Let $X$  be  a  compact set  in $\P^n,$ $n\geq 3.$ 
Assume  $U^+:=\P^n\setminus X$ is pseudoconvex (hence Stein). Then, the following properties hold.
\begin{itemize}
\item[(i)]  Let $\beta$   be  a  smooth $(0,1)$-form on $X$ such that $\dbar \beta=0$ in $\Cc^\infty_{(0,1)}(X).$
Then for each integer $k,$ there is  a  function  $v\in\Cc^k(X)$ such that $\dbar v= \beta.$  
\item[(ii)] If $\mathcal H^2(X)=0,$ then there is a strictly psh function  near $X.$
 \end{itemize}
\end{thmspec}

Recall that on a pseudoconvex domain $U$ in $\P^n,$  Takeuchi \cite{T} and Elencwajg  \cite{E}, proved that if $\delta$  denotes the distance  to the  boundary  of $U$
(with respect to  the Fubini-Study  metric $\omega$), there is a constant $C$  such that near the boundary of $U,$
$$
i\ddbar (-\log \delta)\geq C\omega.
$$

In particular, pseudoconvex  domains are Stein. Indeed, the result is valid when $U$ is pseudoconvex in a compact K\"ahler manifold $M,$ with positive holomorphic bisectional curvature. Moreover the constant $C$ depends only on the curvature. See \cite{E}, in particular, inequality (39), and  Greene-Wu \cite{GW}.

We  will use  the following  result   which is a consequence of Serre's duality  and  H\"ormander's estimates. One  solves the $\dbar$-equation with the  weight $e^\varphi,$ with $\varphi$ psh
instead of the  classical $e^{-\varphi},$ one need a vanishing of the forms on the boundary, see  for example  \cite{BS}.
\renewcommand{\thethmspec}{Theorem 5.2}
\begin{thmspec}
Let $U$  be  a  pseudoconvex domain   in $\P^n.$  Let $\alpha\in L^2_{p,q}(U,\loc),$ $q<n,$  be  a  $\dbar$-closed  form  such that for a given $s >0$
$$
\int |\alpha|^2{1\over \delta^{s+2}} d\lambda <\infty. 
$$
Then  there exists $u\in L^2_{p,q-1}(U,\loc)$ such that
$$
\dbar u=\alpha
$$
and
$$
\int |u|^2{1\over \delta^{s+2}} d\lambda\leq  
{1\over C}\int |\alpha|^2{1\over \delta^{s+2}} d\lambda.
 $$
 Here $d\lambda$  denotes the volume form associated to $\omega.$
\end{thmspec}
\proof
We can consider  that  $\beta$  is extended as a smooth $(0,1)$-form  in $\P^n$
and that $\dbar\beta$  vanishes to infinite order on $X.$  Indeed, the jet of $\beta$ on $X$ satisfies $\dbar\beta=0$ as a jet. So the extension of $\beta$ will satisfy the asserted property. 

Let $\delta$  denote  the  Fubini-Study  distance to $X$ on $U^+.$
We know that $\varphi=-\log\delta$ is psh. We  use  the above H\"ormander's type result. For $s>0,$
there is  an $l\geq 2$ and a form $\psi,$  such that $\dbar \psi=\dbar \beta$ on $U^+,$   with  the  following estimate
$$
\int |\psi|^2{1\over \delta^{s+l}} d\lambda\leq  
C_s\int |\dbar \beta|^2{1\over \delta^{s+l}} d\lambda<\infty.
$$
We can  choose  $s=2n+l.$
Hence for $z\in U^+,$ and $k>0$ fixed,
\begin{eqnarray*}
 {|\psi(z)|^2\over  \delta^k(z)}&\lesssim  & {1\over\delta^{2n+k}(z) }\int_{B(z,\delta(z))}|\psi|^2+   {1\over \delta^k(z)} \sup\limits_{B(z,\delta(z))} |\dbar \psi|\\
 &\lesssim& \int {|\psi(z)|^2\over  \delta^{2n+k}(z)}+o(\delta)\\
 &\lesssim& \int {|\dbar\beta|^2\over  \delta^{2n+l+k}(z)}+o(\delta)\leq C.
\end{eqnarray*}

Hence, $|\psi(z)|=O(\delta^k).$ Consequently,  $\psi$  vanishes on $X$ to any given fixed order. If extended by zero on $X,$ it is  in $\Cc^k(\P^n).$ We can now solve in $\P^n$ the equation
$$
\dbar u=\beta -\psi.
$$
The restriction of $u$ to $X$  satisfies $\dbar u=\beta$  and is  in $\Cc^k(X).$

If  $\mathcal H^2(X)=0,$ then $\omega=\dbar \varphi^{0,1}+\dbar \varphi^{0,1}$ near $X,$ with
$\dbar\varphi^{0,1}=0$ near $X.$ We solve $\dbar u=\varphi^{0,1}$ on $X.$ Then  
$\omega=\ddbar u+\dbar \partial \bar u=i\ddbar \big( {u-\bar u\over i}\big).$  The  function $ v:={u-\bar u\over i}$ is  strictly psh on $X$ and hence in a neighborhood of $X.$
\endproof

\renewcommand{\thethmspec}{Remarks 5.3}
\begin{thmspec}\rm
 \begin{enumerate}
  \item[1.]  A similar  result  can be  obtained for $(p,q)$-forms with $q+1<n.$
  \item[2.] Suppose  $U^+$ is a Stein  domain with smooth boundary and that $H^{2n-2}(U^+)\not=0.$
  As we have seen, there is a strictly psh function on a neighborhood of $\partial U^+.$ Then  a theorem of J. J. Kohn \cite{Ko} asserts that one can solve 
  the $\dbar$-equation in the   Sobolev spaces $H^s(U^+)$ for $s$ large enough. One  can solve it also   in $\Cc^\infty.$
  \item[3.] The  results  are valid  if we replace $\P^n,$ by a  compact simply  connected  K\"ahler manifold $M$, of dimension $n\geq 3,$ with positive holomorphic bisectional curvature, such that $H^{0,1}(M)=0,$  
  and $H^{1,1}(M)$ is  one-dimensional.
\item[4.]  To get  a strictly  psh function  on  a neighborhood of $X,$ it is  enough   to assume the  existence  of a $1$-form $\psi$ near $X$ with
$\partial \psi^{0,1}+\dbar  \psi^{1,0}>0$ on $X$ and  $\dbar \psi^{0,1}=0$ on $X.$
This is satisfied  when $\mathcal H^2(X)=0.$

 Observe however that the existence of such a $1$-form $\psi$
implies the following geometric condition on $X.$  There is no closed current $T$ of dimension $2$, with $\{T\}\not=0,$ such that the component of bi-dimension $(1,1), T_{1,1}$ is positive. Indeed, since $T$ is closed and on $X,$ it follows that 
$\ddbar  (T_{1,1})=0.$ We then have, since $\dbar \psi^{0,1}=0,$
$$
0=\langle T,d\psi\rangle= \langle   T_{1,1},\partial\psi^{0,1} + \dbar \psi^{1,0}\rangle.
$$
This implies that  $T_{1,1}=0,$ and hence $\{T\}=0.$
\end{enumerate}

\end{thmspec}


\section{The case of surfaces}\label{S: dim2}

   Let $U^-$  be a pseudoconvex domain in $\P^2$ with $\Cc^2$-boundary. Let $r$ be a  defining  function for the  boundary  $Z=\partial U^-.$ 
If  for every $z\in Z,$ 
$$\langle i\ddbar r(z),it\wedge\bar{t} \rangle= 0\qquad\text{if}\qquad  \langle\partial r(z),t  \rangle=0$$
then we say that the  boundary is Levi flat.

It is  not known if such domains  exist, even if we  assume that the  boundary is  real analytic.
However for arbitrary K\"ahler  surfaces $M$ and
for a Levi  flat surface, there is a positive current, $\ddbar$-closed of mass 1 directed by the foliation on the  boundary $Z.$
 It is    shown in \cite{FS}  for $\P^2$ and in \cite{DNS} in general, that when there is no  positive closed  current of mass 1  directed  by the  foliation, then $T$ is  unic.
As we have seen in  Proposition 2.2, for $\P^2$  or more generally  for simply connected K\"ahler surfaces, there is no positive closed current on $Z.$ Indeed, we can assume that $\overline {U^-}$ satisfies $H_{dR}^2(\overline{ U^-})=0.$

\renewcommand{\thethmspec}{Theorem 6.1}
\begin{thmspec}
Let  $X$  be a compact  set in a compact K\"ahler surface $M$, such that $ \mathcal H^2(X)=0.$ Assume $X$ supports a positive $\ddbar$-closed  current $T$ of mass $1.$
Then there is a real analytic $(0,1)$-form $\varphi^{0,1},$ $\dbar$-closed in a neighborhood of $X$ and  such that there is  no function $u$
in the  Sobolev space  $W{^2}(M),$  with $\dbar u=\varphi^{0,1}$ on $X.$

In particular if $X$ is a real hypersurface, there is no solution $u \in W{^ \frac {3} {2}}(X),$ for the equation $\dbar u=\varphi^{0,1}.$ 
\end{thmspec}

\proof
For arbitrary $X$ the meaning of $\dbar u=\varphi^{0,1}$ on $X,$ is that for every current $R \in W^{-1}(M),$ supported on $X,$ 
$\langle R,\varphi^{0,1}\rangle= \langle  R, \dbar u\rangle.$ If $X$ is a $\Cc^1,$ hypersurface, a function in $u \in W{^ \frac {3} {2}}(X) ,$ extends as a function in  $W{^2}(M),$  and a function in $W{^2}(M),$ restricts to  $W{^ \frac {3} {2}}(X),$ so  $\dbar u=\varphi^{0,1}$ on $X,$ makes sense and the two notions coincide.

Since $\mathcal H^2(X)=0,$ the current $T$ is not closed and hence $\partial T$ is non-zero.
It is   shown  in \cite{FS} that if $T\geq 0,$ $\ddbar T=0$ and  $\int T\wedge \Omega=1$  then  $T=\Omega+\partial S+\overline{\partial S}+i\ddbar v,$
with $S,$ $\partial S,$ $\dbar S\in L^2$ and $v\in L^p$ for all $p<2.$  Here $\Omega$ is smooth and represents the class of $T.$

It follows that  
\[  \partial T=-\dbar\partial \bar S\]
is in the Sobolev  space $W^{-1}(M).$  If $u\in W^2(M),$ since $\ddbar T=0,$  
\[ \langle \partial T, \dbar u\rangle=0.\]
On the other hand in a neighborhood of $X,$  $\omega=\partial \varphi^{0,1}+\dbar \varphi^{1,0},$ $\dbar \varphi^{0,1}=0.$ Hence
$$
1=\langle T,\omega\rangle=-2\Re \langle \partial T,\varphi^{0,1}\rangle.
$$
 So $\langle \partial T,\varphi^{0,1}\rangle\not=0.$ Hence we cannot have $\varphi^{0,1} = \dbar u$ on $X,$ 
  with $u$ having an extension in $W{^2}(M).$ 
 \endproof

\renewcommand{\thethmspec}{Remark 6.2}
\begin{thmspec} \rm

\item If $X=\partial U^-,$ is Levi-flat with  $H_{dR}^2(\overline{ U^-})=0,$ then  $\partial T$  is  of order  zero \cite {FS} , and the  same proof  shows there is  no  continuous function $u$ on $\overline{U^-}$
 such that  $\dbar  u=\varphi^{0,1}.$

\end{thmspec}

\renewcommand{\thethmspec}{Question}\rm 
\begin{thmspec} 
Suppose   $U$  is  a  smooth pseudoconvex   domain  in $\P^n,$ $n\geq 2.$  
Assume  
$H_{dR}^2( {U})=0.$  Is there   a   strictly psh function  near the boundary?
\end{thmspec}
There is an example of a compact K\"ahler  surface $M,$ with a Stein domain $U \subset M$
with real analytic boundary, but all bounded psh functions in $U$ are constant,\cite {OS} .

\renewcommand{\thethmspec}{Theorem 6.3}
\begin{thmspec}
Let $M$ be a compact complex manifold of dimension $n$. Let $Z$ be an irreducible real-analytic set of real dimension $m.$ Suppose that there is in $Z$ a germ of a complex analytic set of dimension 
$p\geq 1.$ Define $A_p$ as the set of points $z\in Z,$ with a germ of complex analytic set of dimension $\geq p$ through $z.$ 
 Then $A_p$ is closed.
 If there is a strictly psh function near $Z$ then $A_1$ is empty.
\end{thmspec}
\proof
It is possible to cover $Z$ with finitely many open sets $(V_i)$ such that on $V_i,$ the equation of $Z,$ is
 $\rho_i (z,\overline z)=0.$  Here $\rho_i(z,\overline z)$ a real-valued analytic function in $z,\overline z$. We can assume that the functions  $\rho_i(z,w)$  are holomorphic in 
 $V_i\times V_i'$ , where $V_i'$ is the image of $V_i$ by conjugation.
 
 Consider a germ of a complex analytic set $L$ parametrized by a holomorphic map $u:\D^p\to L$. Then, the function $\rho_i(u(t),\overline u(s))$ on $\D^{2p}$ vanishes when $s=\overline t$, where $\overline u(s):=\overline{u(\overline s)}$. Since this function is holomorphic, its zero set is a complex analytic set. Hence, it vanishes everywhere in $\D^{2p}$. Fixing an arbitrary $s$, we deduce that
$\rho_i(z, u(s))=0$ for $z\in L$. 

Define $ L' := \left\lbrace z\in V_i, \bigcap_{s} \rho_i(z, u(s))=0 \right\rbrace. $
Then $L'$ is an extension of the germ $L$ to an analytic set in $V_i,$ 
contained in $Z.$ The size of the $V_i$'s is fixed. It follows easily that $A_p$ is closed. This is precisely the Segre argument, see \cite { DS2} and \cite[Example 7]{FS2}. 

If $A_1$ is non-empty, then it satisfies the local maximum principle for psh, functions in a neighborhood. If $u$ is strictly psh in a neighborhood of $Z,$ it reaches it's maximum on 
$A_1$ at a point $p.$ As in Remark 3.4 (4), we can arrange that the maximum is strict. A contradiction.
\endproof
\renewcommand{\thethmspec}{Theorem 6.4}
\begin{thmspec}

Let $M$ be a compact complex surface. Let $U$ be a smooth  domain with connected real analytic boundary $Z.$ Assume that $Z$ admits a point of strict pseudoconvexity. Then either there is a compact complex curve on the boundary, or there is  a   strictly psh function  near the boundary. In particular if $Z \subset \P^2,$
there is  a   strictly psh function in a neighborhood of $Z.$
\end{thmspec}
\proof
Suppose there is no compact complex curve on the boundary.
Let  $S$ be  the union of strictly pseudoconvex points and strictly
pseudoconcave points on $Z.$  Let $W= Z\setminus S$. Since the defining function $r$ is real analytic, then $W$ is a real analytic set of dimension $\leq2$.  It admits a stratification by smooth manifolds. If there is a curve of real dimension $2$
with complex tangents, then it is a complex curve. By the above theorem it has no boundary and it is necessarily of finite area, since this is the case for $W.$ Hence it is a compact complex curve. Since this is not possible, then the set $C$ with complex tangents is at most of finite one dimensional Hausdorff measure. 

On the other hand, if there is no strictly psh function near the boundary, there is  a  positive $\ddbar$-closed  current $T$ of bi-dimension $(1,1)$ of mass 1 supported on $\partial U.$ Such a current has no mass on the set of strictly pseudoconvex points, nor on the set of strictly pseudoconcave points, as follows from equations (3.1). Hence it  is supported on $W.$ Since $T$ is of bi-dimension$(1,1)$ it is supported on $C$. But such currents don't give mass to sets of  2-Hausdorff 
dimension zero.  See \cite {S3} for more details on the geometry of such currents.  It follows that there is a strictly psh function near $\partial U.$
\endproof
\renewcommand{\thethmspec}{Remark}
\begin{thmspec} \rm
It is shown in \cite { DS2} and \cite[Example 7]{FS2}  that no Levi-flat  hypersurface $Z$ exists in $\P^2$ if we assume it is smooth and real algebraic. Hence, there is a point of strict pseudoconvexity. It follows that for smooth real algebraic surfaces in $\P^2$
there is a strictly psh function in a neighborhood of $Z.$

\end{thmspec}

\renewcommand{\thethmspec}{Theorem 6.5}
\begin{thmspec}
Let  $X$  be a compact  set in a compact K\"ahler surface $M$.
\begin{enumerate}
 \item [1.]  If $X$ is (locally) pluripolar, then either there is a positive closed current $T$ of bi-degree $(1,1),$ and of mass one supported on $X,$ or there is a strictly psh function in a neighborhood of $X.$
 
 \item[2.] If $X$ is of  Lebesgue measure zero, either there is a strictly psh function in a neighborhood of $X$ or a positive closed current $T$ of mass one supported on $X,$ or there is a $(2,0)$-form $A$ in $L_{(2,0)}^2 (M)$ of norm $1$ and holomorphic in $M\setminus X$.
 \end{enumerate}

\end{thmspec}
\proof
Assume $X$ is pluripolar. If there is no strictly psh function near $X,$ then there is a positive $\ddbar$-closed current $T$ of mass $1,$ supported on $X.$ 

Moreover, as we have seen that,
 $T=\Omega+\partial S+\overline{\partial S}+i\ddbar u,$
with $S,$ $\partial S,$ $\dbar S\in L^2$ and $u\in L^p$ for all $p<2.$ Hence $\partial T=-\dbar\partial \bar S.$ So $\partial \bar S$ is holomorphic out of $X$, which is pluripolar.
Since holomorphic functions  in $L^2,$ extend holomorphically through pluripolar sets, it follows that $\dbar\partial \bar S=0$ in $M$. Hence $T$ is closed.

If we assume only that $X$ is of Lebesgue measure zero, and there is no positive closed current supported on $X$, then $A=\partial \bar S$ is a non-identically zero, holomorphic form in $M\setminus X$, which is in $L^2$.
\endproof


\end{document}